\documentclass[14pt,a4paper]{extreport}

\usepackage{amsfonts,amsmath,amsthm,amscd,amssymb,latexsym,cite}

    \usepackage[T2A]{fontenc}
    \usepackage[cp1251]{inputenc}
    \usepackage[russian]{babel}
\usepackage[matrix, arrow, curve]{xy}
\usepackage{amsfonts}
\usepackage{amssymb}
\usepackage{stmaryrd}
\usepackage{textcomp}
\usepackage{amsthm}
\usepackage{amsmath}

\usepackage{ifpdf}
\ifpdf
\usepackage[pdftex]{graphicx}
\usepackage[pdftex,dvipsnames,usenames]{color}
\else
\usepackage{graphicx}
\usepackage[dvips,dvipsnames,usenames]{color}
\fi

\newtheorem{theorem}{Теорема}
\newtheorem{lemma}{Лемма}

\newtheorem{definition}{Определение}
\newtheorem{corollary}{Следствие}

\newtheorem{remark}{Замечание}

\makeatletter 
\def\@biblabel#1{#1. }
\makeatother

\begin{document}

\subsection*{\center Pretangent spaces with nonpositive and nonnegative Aleksandrov curvature}\begin{center}\textbf{\textmd{\textbf{V. Bilet and O. Dovgoshey}}} \end{center}
\parshape=5
1cm 13.5cm 1cm 13.5cm 1cm 13.5cm 1cm 13.5cm 1cm 13.5cm \noindent \small {\bf Abstract.}
  We find conditions under which the pretangent spaces to general metric spaces have the nonpositive Aleksandrov curvature or nonnegative one. The infinitesimal structure of general metric cpaces with Busemann convex pretangent spaces is also described.

\bigskip

\parshape=1
0.3 cm 13.5cm \textbf{2010 MSC: 54E35.}

\medskip

\parshape=2
1cm 13.5cm 1cm 13.5cm  \noindent \small {\bf Key words:} pretangent space, CAT(0)-space, Aleksandrov curvature, Busemann convexity, infinitesimal geometry of metric spaces.

\subsection*{\center Предкасательные пространства с неположительной и неотрицательной по Александрову кривизной}\begin{center}\textbf{\textmd{\textbf{\normalsize В. В. Билет и А. А. Довгошей}}}\end{center}
\parshape=5
1cm 13.5 cm 1cm 13.5cm  1cm 13.5 cm  1cm 13.5 cm  1cm 13.5 cm  \noindent \small {\bf Аннотация.}
 Мы находим условия, при которых предкасательные пространства к общим метрическим пространствам имеют неположительную или неотрицатель\-ную кривизну по Александрову. Также описана инфинитезимальная структура общих метрических пространств с выпуклыми по Буземанну предкасательными пространствами.

\bigskip

\parshape=1
0.3 cm 13.5 cm \textbf{2010 MSC: 54E35.}

\medskip

\parshape=2
1cm 13.5cm 1cm 13.5cm  \noindent \small {\bf Ключевые слова:} предкасательное пространство, CAT(0)-пространство, кривизна по Александрову, выпуклость по Буземанну, инфинитезимальная геометрия метрических пространств.

\newpage \noindent{\bf 1. Введение.}
Предкасательные и касательные пространства, используемые в настоящей работе, были введены в \cite{DM} (см. также \cite{DM1}). Напомним необходимые определения. 

Пусть $(X,d,p)$ $-$ метрическое пространство с отмеченной точкой $p.$ Зафиксируем
последовательность $\tilde r$ положительных вещественных чисел $r_{n},$
стремящихся к нулю. Назовём $\tilde r$ \emph{нормирующей}
последовательностью. Будем обозначать через $\tilde X$ множество всех
последовательностей точек из $X.$

\begin{definition}  Две последовательности $\tilde{x}, \tilde{y} \in \tilde{X},$
$\tilde{x}=\{x_n\}_{n\in \mathbb N}$ и $\tilde{y}=\{y_n\}_{n\in \mathbb N},$
взаимностабильны относительно нормирующей последовательности
$\tilde{r}=\{r_n\}_{n\in \mathbb N},$ если существует конечный предел
\begin{equation}\label{eq1}
\lim_{n\to\infty}\frac{d(x_n,y_n)}{r_n}:=\tilde{d}_{\tilde r}(\tilde{x},\tilde{y})=\tilde{d}(\tilde{x},\tilde{y}).
\end{equation}
\end{definition}

Cемейство $\tilde{F}\subseteq \tilde{X}$ \emph{самостабильное}, если любые две
последовательности $\tilde{x}, \tilde{y} \in\tilde{F}$ взаимностабильны,
$\tilde{F}\subseteq \tilde{X}$ $-$ \emph{максимальное самостабильное}, если $\tilde{F}$
самостабильное и для произвольной $\tilde{z} \in \tilde{X}\setminus \tilde{F}$
существует $\tilde{x}\in\tilde{F}$ такая, что $\tilde{x}$ и $\tilde{z}$ не
взаимностабильны. Из леммы Цорна легко следует, что для каждой нормирующей
последовательности $\tilde{r}=~\{r_n\}_{n\in \mathbb N}$ существует максимальное
самостабильное семейство $\tilde X_{p,\tilde{r}}$ такое, что постоянная
последовательность $\tilde p=\{p,p,...\}\in {\tilde X_{p,\tilde{r}}}.$

Рассмотрим функцию $\tilde d: \tilde X_{p, \tilde r}\times \tilde X_{p, \tilde r}\to
\mathbb R,$ где $\tilde d(\tilde x, \tilde y)=\tilde d_{\tilde r}(\tilde x, \tilde y)$
определена через $\eqref{eq1}.$ Очевидно, $\tilde {d}$ симметрична и неотрицательна.
Кроме того, из неравенства треугольника для $d$ имеем $\tilde d (\tilde x, \tilde y)
\le \tilde d (\tilde x, \tilde z) + \tilde d(\tilde z, \tilde y)$ для всех $\tilde x,
\tilde y, \tilde z$ из $\tilde X_{p,\tilde{r}}.$ Следовательно $(\tilde X_{p,\tilde{r}},
\tilde d)$ $-$ псевдометрическое пространство.

Определим отношение эквивалентности $\sim$ на $\tilde X_{p,\tilde{r}}$ как $ \tilde x
\sim \tilde y $  тогда и только тогда, когда $\tilde d_{\tilde r}(\tilde x, \tilde
y)=0.$ Обозначим через $\Omega_{p,\tilde r}^{X}$ множество всех
классов эквивалентности на $\tilde X_{p,\tilde{r}},$ порождённых отношением $\sim.$ Для
$\alpha, \beta \in {\Omega_{p,\tilde r}^{X}}$ положим $\rho(\alpha,\beta)=\tilde
d(\tilde x, \tilde y),$ где $\tilde x \in \alpha$ и $\tilde y \in \beta,$ тогда $\rho$
$-$ метрика на $\Omega_{p,\tilde r}^{X}.$ Переход от псевдометрического пространства
$(\tilde X_{p,\tilde{r}},\tilde d)$ к метрическому пространству $(\Omega_{p,\tilde r}^{X},
\rho)$ будем называть \emph{метрической идентификацией} $(\tilde X_{p, \tilde{r}},\tilde d).$

\begin{definition} Пространство $(\Omega_{p,\tilde r}^{X}, \rho)$ называется предкасательным к $X$ в
точке $p$ относительно нормирующей последовательности $\tilde r.$\end{definition}

Пусть $\{n_{k}\}_{k\in\mathbb N}$ $-$ бесконечная, строго возрастающая последовательность натуральных чисел. Обозначим через $\tilde r'$ подпоследовательность $\{r_{n_k}\}_{k\in\mathbb N}$ нормирующей последовательности $\tilde r=\{r_n\}_{n\in\mathbb N}$ и пусть $\tilde x':=\{x_{n_{k}}\}_{k\in\mathbb N}$ для каждой $\tilde x=\{x_n\}_{n\in\mathbb N}\in \tilde X.$ Ясно, что если $\tilde x$ и $\tilde y$ взаимностабильны относительно $\tilde r,$ то $\tilde x'$ и $\tilde y'$ взаимостабильны относительно $\tilde r'$ и
$\tilde d_{\tilde r} (\tilde x, \tilde y) = \tilde d_{\tilde r'} (\tilde x', \tilde y').$
Если $\tilde X_{p, \tilde r}$ $-$ максимальное самостабильное относительно $\tilde r$ семейство, тогда, по лемме Цорна, существует максимальное самостабильное относительно $\tilde r'$ семейство $\tilde X_{p, \tilde r'}$ такое, что
$\{\tilde x' : \tilde x \in \tilde X_{p, \tilde r}\}\subseteq \tilde X_{p, \tilde r'}.$

Обозначим через $in_{\tilde r'}$ отображение из $\tilde X_{p, \tilde r}$ в $\tilde X_{p, \tilde r'}$ с $in_{\tilde r'} (\tilde x) = \tilde x'$ для всех $\tilde x \in \tilde X_{p, \tilde r}.$ После метрической идентификации отображение $in_{\tilde r'}$ переходит в изометрическое вложение $em' : \Omega_{p, \tilde r}^{X} \rightarrow \Omega_{p, \tilde r'}^{X},$ для которого диаграмма
\begin{equation}\label{diag}
\begin{array}{ccc}
\tilde X_{p, \tilde r} & \xrightarrow{\ \ \mbox{\emph{in}}_{\tilde r'}\ \ } &
\tilde X_{p, \tilde r^{\prime}} \\
\!\! \!\! \!\! \!\! \! \pi\Bigg\downarrow &  & \! \!\Bigg\downarrow \pi^{\prime}
\\
\Omega_{p, \tilde r}^{X} & \xrightarrow{\ \ \mbox{\emph{em}}'\ \ \ } & \Omega_{p, \tilde
r^{\prime}}^{X}
\end{array}
\end{equation} коммутативна. Здесь $\pi, \pi'$ отображения проектирования на соответствующие факторпространства, $\pi (\tilde x):=\{\tilde y\in\tilde X_{p,\tilde r} : \tilde d_{\tilde r}(\tilde x, \tilde y)=0\}$ и $\pi' (\tilde x):=\{\tilde y\in\tilde X_{p,\tilde r'} : \tilde d_{\tilde r'}(\tilde x, \tilde y)=0\}.$

Предкасательное пространство $\Omega_{p,\tilde r}^{X}$ является \emph{касательным}, если $em':\Omega_{p,\tilde r}^{X}\rightarrow \Omega_{p,\tilde r'}^{X}$ биективно для каждого $\tilde X_{p,\tilde r'}.$

\noindent{\bf 2. Постановка задачи.} В настоящей работе исследуются условия на метрическое пространство $X,$ при которых предкасательные пространства $\Omega_{p,\tilde r}^{X}$ имеют неотрицательную по Александрову кривизну и условия, при которых кривизна по Александрову пространств $\Omega_{p,\tilde r}^{X}$ неположительна, т.е. $\Omega_{p,\tilde r}^{X}$ являются CAT(0)-\emph{пространствами.} Хорошо известные определения неотрицательности (неположительности) кривизны по Александрову даются для геодезических пространств с использованием так называемых треугольников сравнения (comparison triangles) см., например, в \cite[гл. 4]{BBI} и являются достаточно громоздкими. Поэтому прямое использование этих определений при исследовании предкасательных пространств к общим метрическим пространствам представляется затруднительным. Положение изменилось после того, как Берг и Николаев дали следующую характеристику CAT(0)-пространств.

\begin{theorem}\emph{\cite{BN}}
Пусть $(X,d)$ $-$ геодезическое пространство. $X$ является CAT(0)-пространством тогда и только тогда, когда неравенство четырехугольника
\begin{equation}\label{eq3}
d^{2}(w,y)+d^{2}(x,z)\le d^{2}(w,x)+d^{2}(x,y)+d^{2}(y,z)+d^{2}(z,w)
\end{equation}
выполнено для любых точек $w,x,y,z\in X.$
\end{theorem}

Простое доказательство этой теоремы было найдено в \cite{Sato}.
После появления \cite{BN}, Лебедева и Петрунин \cite{LP} получили аналогичную характеристику пространств с неотрицательной по Александрову кривизной.

\begin{theorem} \emph{\cite{LP}}
Пусть $(X,d)$ $-$ полное геодезическое пространство. $X$ является пространством неотрицательной по Александрову кривизны тогда и только тогда, когда неравенство
\begin{equation}\label{eq3*}
\frac{1}{3}(d^{2}(x,y)+d^{2}(y,z)+d^{2}(z,x))\le d^{2}(w,x)+d^{2}(w,y)+d^{2}(w,z)
\end{equation}
выполнено для любых точек $w,x,y,z\in X.$
\end{theorem}

Таким образом, для описания структуры метрических пространств $(X,d),$ предкасательные к которым принадлежат классу CAT(0) или имеют неотрицательную по Александрову кривизну, достаточно найти:
\newline(i) инфинитезимальные аналоги неравенств \eqref{eq3} и \eqref{eq3*}; \newline(ii) условия геодезичности $\Omega_{p,\tilde r}^{X};$ \newline(iii) условия полноты $\Omega_{p,\tilde r}^{X}.$

Пункт (i) можно реализовать достаточно просто, если использовать так называемый ``принцип переноса'', доказанный в \cite{BD}.

\textbf{Принцип переноса и его применение.}

Пусть $(X,d)$ -- метрическое пространство. Для каждого $n\in\mathbb
N$ обозначим через $X^{n}$ множество всех $n$-наборов $x=(
x_1,...,x_n)$ таких, что $x_{i}~\in~X,$ $i=~1,...,n.$
Обозначим через $\textbf{\emph{M}}_{n}$
пространство вещественных $(n~\times~n)$-матриц $\textbf{\emph{t}}$ с топологией
поточечной сходимости. Пусть $\mathfrak{M}$ -- фиксированный класс непустых
метрических пространств и пусть $\mathfrak{F}$ -- фиксированное семейство
непрерывных однородных функций $f: \textbf{\emph{M}}_{n}\rightarrow
\mathbb R, n=n(f)$ степени однородности $s=s(f)>0.$

Будем говорить, что $\mathfrak{M}$ определяется семейством
$\mathfrak{F},$ если следующие два условия эквивалентны для всякого
метрического пространства $(X,d):$

$\bullet$ $(X, d)\in\mathfrak{M};$

$\bullet$ неравенство $f(\textbf{\emph{m}}(x_1,...,x_n))\ge 0$ выполнено для
каждой $f\in\mathfrak{F}$ и всех
\begin{equation*}
\textbf{\emph{m}}(x_1,
x_2,...,x_n)=\left(\begin{array}{cccc}
   d(x_{1}, x_{1})&d(x_{1}, x_{2})&...&d(x_{1}, x_{n})\\
   d(x_{2}, x_{1})&d(x_{2}, x_{2})&...&d(x_{2}, x_{n})\\
    \vdots&\vdots&\ddots&\vdots\\
   d(x_{n}, x_{1})&d(x_{n}, x_{2})&...&d(x_{n}, x_{n})\\
    \end{array}\right), (x_1, x_2, ..., x_n)\in X^{n}.
\end{equation*}

Пусть $(X,d, p)$ -- метрическое пространство c отмеченной точкой $p.$  Полагаем
\begin{equation}\label{norm}
\delta_{p}(x_1,...,x_n):=\max_{1\le i\le n}d(x_{i},p)
\end{equation} для $(x_1,...,x_n)\in X^{n}.$ Для $f\in\mathfrak M$ определим функцию $f^{*}: X^{n}\rightarrow \mathbb R$ правилом
\begin{equation*}\label{Func}
f^{*}(x_1,x_2,...,x_n):=\begin{cases}
         f\left(\frac{\textbf{\emph{m}}(x_1, x_2,...,x_n)}{\delta_{p}(x_1, x_2,...,x_n)} \right), & \mbox{если} $ $ (x_1,x_2,...,x_n) \ne (p,p,...,p)\\
         0,& \mbox{если}$ $ (x_1,x_2,...,x_n) = (p,p,...,p). \\
         \end{cases}
\end{equation*}

\begin{lemma} \emph{\cite{BD}}
Пусть $(X,d,p)$ -- метрическое пространство с отмеченной точкой $p$ и
пусть $\mathfrak{M}$ -- семейство метрических пространств,
определяемое некоторым семейством $\mathfrak{F}.$ Следующие два утверждения
эквивалентны:\newline \emph{(i)} Каждое предкасательное пространство
$\Omega_{p,\tilde r}^{X}$ принадлежит $\mathfrak{M};$ \newline \emph{(ii)}
Неравенство \begin{equation*}\label{eq1.10} \liminf_{x_1,x_2,...,x_n
\to p}f^{*}(x_1,x_2,...,x_n)\ge 0 \end{equation*} выполняется для
всякой $f\in \mathfrak{F}.$
\end{lemma}

Положим

\begin{equation}\label{eq5}
A_{1}(w,x,y,z):= \frac{d^{2}(w,x)+d^{2}(x,y)+d^{2}(y,z)+d^{2}(z,w)-(d^{2}(w,y)+d^{2}(x,z))}{(\delta_{p}(w,x,y,z))^{2}},
\end{equation}
\begin{equation}\label{eq6}
A_{2}(w,x,y,z):= \frac{d^{2}(w,x)+d^{2}(w,y)+d^{2}(w,z)-\frac{1}{3}(d^{2}(x,y)+d^{2}(y,z)+d^{2}(z,x))}{(\delta_{p}(w,x,y,z))^{2}}
\end{equation}
при $(w,x,y,z)\ne (p,p,p,p)$ и
\begin{equation*}
A_{1}(p,p,p,p)=A_{2}(p,p,p,p)= 0.
\end{equation*}

Следующая лемма следует непосредственно из леммы 1 при $f^{*}=A_{1}$ и $f^{*}=A_{2}.$

\begin{lemma}
Пусть $(X,d,p)$ $-$ метрическое пространство с отмеченной точкой $p.$ Неравенство
\begin{equation*}
\rho^{2}(\alpha,\beta)+\rho^{2}(\gamma,\eta)\le\rho^{2}(\alpha,\gamma)+\rho^{2}(\alpha,\eta)+\rho^{2}(\beta,\gamma)+\rho^{2}(\beta,\eta)
\end{equation*}
выполнено для любого предкасательного пространства $\Omega_{p,\tilde r}^{X}$ и любых $\alpha, \beta, \gamma, \eta\in\Omega_{p,\tilde r}^{X}$ тогда и только тогда, когда
\begin{equation*}\liminf_{w,x,y,z\to p}A_{1}(w,x,y,z)\ge 0.
\end{equation*}
Неравенство
\begin{equation*}
\frac{1}{3}(\rho^{2}(\alpha,\beta)+\rho^{2}(\beta,\gamma)+\rho^{2}(\gamma,\alpha))\le\rho^{2}(\alpha,\eta)+\rho^{2}(\beta,\eta)+\rho^{2}(\gamma,\eta)
\end{equation*} выполнено для любого предкасательного пространства $\Omega_{p,\tilde r}^{X}$ и любых $\alpha, \beta, \gamma, \eta\in\Omega_{p,\tilde r}^{X}$ тогда и только тогда, когда
\begin{equation*}\label{eq4*}\liminf_{w,x,y,z\to p}A_{2}(w,x,y,z)\ge 0.
\end{equation*}
\end{lemma}

\textbf{Геодезичность предкасательных пространств.}
Условия геодезичности предкасательных пространств к общему метрическому пространству исследовались в работе \cite{Bilet}. Для удобства читателя приведем необходимые для дальнейшего результаты этой работы.

Пусть $(X,d,p)$ $-$ метрическое пространство с отмеченной точкой $p.$ Напомним, что точка $m\in X$ называется \emph{срединной точкой} для $x,y\in X,$ если $d(m,x)=d(m,y)=$ $=\frac{1}{2}d(x,y).$

\begin{definition}
Будем говорить, что пространство $X$ является срединно выпуклым в точке $p,$ если для любых двух последовательностей $\tilde x=\{x_n\}_{n\in\mathbb N}$ и $\tilde y=\{y_n\}_{n\in\mathbb N}\in\tilde X,$ сходящихся к $p,$ найдется последовательность $\tilde z=\{z_n\}_{n\in\mathbb N}\in\tilde X,$ $\tilde z=\tilde z(\tilde x, \tilde y)$ такая, что
\begin{equation}\label{eq7}
d(x_n, z_n)=\frac{1}{2}d(x_n,y_n)+o(\delta_{p}(x_n,y_n))
\quad \mbox{и} \quad
d(y_n, z_n)=\frac{1}{2}d(x_n,y_n)+o(\delta_{p}(x_n,y_n)),
\end{equation}
где $\delta_{p}(x,y)$ определена формулой \eqref{norm}, а формулы \eqref{eq7} означают, что
\begin{equation}\label{eq9}
\lim_{n\to\infty}\frac{\left|d(x_n,z_n)-\frac{1}{2}d(x_n,y_n)\right|}{\delta_{p}(x_n,y_n)}=\lim_{n\to\infty}\frac{\left|d(y_n,z_n)-\frac{1}{2}d(x_n,y_n)\right|}{\delta_{p}(x_n,y_n)}=0.
\end{equation}
\end{definition}

Последовательность $\tilde z$ в определении 3 будем называть \emph{инфинитезимальной срединной точкой} для $\tilde x$ и $\tilde y.$

\begin{remark}
При $d(x_n,p)=d(y_n,p)=0$ формула \eqref{eq9} не определена, но ее легко доопределить, считая, что выражения под знаком пределов равны нулю при $x_n=$ $=y_n=z_n=p$ и $+\infty$ при $x_n=y_n=p, z_n\ne p.$ Таким образом, если $\tilde x=\tilde y=\tilde p,$ где $\tilde p~=~(p,p,p,...)$ $-$ стационарная последовательность, а $\tilde z=\{z_n\}_{n\in\mathbb N}$ $-$ инфинитезимальная срединная точка для $\tilde x$ и $\tilde y,$ то $z_n=p$ для всех достаточно больших $n.$
\end{remark}

\begin{lemma} \emph{\cite{Bilet}}
Пусть $(X,d,p)$ $-$ метрическое пространство с отмеченной точкой $p.$ Если $X$ является срединно выпуклым в точке $p,$ то любое сепарабельное касательное пространство $\Omega_{p,\tilde r}^{X}$ является геодезическим.
\end{lemma}

\begin{lemma} \emph{\cite{Bilet}}
Пусть $(X,d,p)$ $-$ метрическое пространство с отмеченной точкой $p.$ Если все предкасательные пространства $\Omega_{p,\tilde r}^{X}$ являются геодезическими, то $X$ $-$ срединно выпукло в точке $p.$
\end{lemma}

\textbf{Полнота предкасательных пространств.} На сегодняшний день нет ни одного примера метрического пространства $X$ с предкасательным пространством, не являющимся полным, однако полнота произвольного предкасательного пространства остается не доказанной. В случае касательных $\Omega_{p,\tilde r}^{X}$ справедлива следующая

\begin{lemma} \emph{\cite{DM}}
Пусть $(X,d,p)$ $-$ метрическое пространство с отмеченной точкой $p.$ Тогда любое касательное пространство $\Omega_{p,\tilde r}^{X}$ является полным.
\end{lemma}

\noindent{\bf 3. Инфинитезимальные версии теорем Берга-Николаева и Лебедевой-Петру\-нина}

\begin{theorem}
Пусть $(X,d,p)$ $-$ метрическое пространство с отмеченной точкой $p$. Если $X$ является срединно выпуклым в точке $p$ и
\begin{equation}\label{eq10}\liminf_{w,x,y,z\to p}A_{1}(w,x,y,z)\ge 0,
\end{equation} где функция $A_{1}(w,x,y,z)$ определена формулой \eqref{eq5}, то любое сепарабельное касательное пространство $\Omega_{p,\tilde r}^{X}$ является CAT(0)-пространством.
\end{theorem}

\begin{proof}
Пусть $X$ $-$ срединно выпукло в точке $p$ и имеет место \eqref{eq10}. Рассмотрим произвольное сепарабельное касательное пространство $\Omega_{p,\tilde r}^{X}$ с метрикой $\rho.$ По лемме 3 пространство $\Omega_{p,\tilde r}^{X}$ является геодезическим, а по лемме 2 неравенство
\begin{equation}\label{eq9*}
\rho^{2}(\alpha,\beta)+\rho^{2}(\gamma,\eta)\le\rho^{2}(\alpha,\gamma)+\rho^{2}(\alpha,\eta)+\rho^{2}(\beta,\gamma)+\rho^{2}(\beta,\eta)
\end{equation} имеет место для любых $\alpha, \beta, \gamma, \eta\in\Omega_{p,\tilde r}^{X}.$ Следовательно, по теореме 1, $\Omega_{p,\tilde r}^{X}$ является CAT(0)-пространством.
\end{proof}

\begin{theorem}
Пусть $(X,d,p)$ $-$ метрическое пространство с отмеченной точкой $p$. Если все предкасательные пространства $\Omega_{p,\tilde r}^{X}$ являются CAT(0)-пространствами, то $X$ $-$ срединно выпукло в точке $p$ и имеет место соотношение \eqref{eq10}.
\end{theorem}

\begin{proof}
Пусть все предкасательные пространства являются CAT(0)-прост\-ранствами. По определению, любое CAT(0)-пространство является геодезическим. Следовательно, все $\Omega_{p,\tilde r}^{X}$ являются геодезическими. По лемме 4, $X$ $-$ срединно выпукло в точке $p.$ По теореме 1 из принадлежности $\Omega_{p,\tilde r}^{X}\in CAT(0)$ следует неравенство \eqref{eq9*} для любых $\alpha, \beta, \gamma, \eta\in\Omega_{p,\tilde r}^{X}.$ Значит \eqref{eq9*} имеет место для всех $\Omega_{p,\tilde r}^{X}$ и всех $\alpha, \beta, \gamma, \eta\in\Omega_{p,\tilde r}^{X}.$ Последнее, по лемме 2, равносильно выполнению \eqref{eq10}.
\end{proof}

\begin{theorem}
Пусть $(X,d,p)$ $-$ метрическое пространство с отмеченной точкой $p$. Если $X$ является срединно выпуклым в точке $p$ и
\begin{equation}\label{eq10*}\liminf_{w,x,y,z\to p}A_{2}(w,x,y,z)\ge 0,
\end{equation} где функция $A_{2}(w,x,y,z)$ определена формулой \eqref{eq6}, то любое сепарабельное касательное пространство $\Omega_{p,\tilde r}^{X}$ является пространством неотрицательной по Александрову кривизны.
\end{theorem}

Доказательство этой теоремы аналогично доказательству теоремы 3. Заметим только, что вместо теоермы 1 (Берга-Николаева) для произвольных геодезических пространств нужно использовать теорему 2 (Лебедевой-Петрунина), в которой предполагается полнота рассматриваемого геодезического пространства. Использование теоремы 2 возможно, так как лемма 5 гарантирует полноту касательных $\Omega_{p,\tilde r}^{X}$.

Следующая теорема полностью аналогична теореме 4.

\begin{theorem}
Пусть $(X,d,p)$ $-$ метрическое пространство с отмеченной точкой $p$. Если все предкасательные пространства $\Omega_{p,\tilde r}^{X}$ являются полными геодезическими пространствами неотрицательной по Александрову кривизны, то $X$ $-$ срединно выпукло в точке $p$ и имеет место соотношение \eqref{eq10*}.
\end{theorem}

Из теорем 3 - 6 получаем следующее

\begin{corollary}
Пусть $(X,d,p)$ $-$ метрическое пространство с отмеченной точкой $p$. Предположим, что любое предкасательное пространство $\Omega_{p,\tilde r}^{X}$ является сепарабельным и касательным. В этом случае:

(i) $X$ $-$ срединно выпукло в точке $p$ и выполнено \eqref{eq10} тогда и только тогда, когда все $\Omega_{p,\tilde r}^{X}$ $-$ геодезические пространства неположительной по Александрову кривизны.

(ii) $X$ $-$ срединно выпукло в точке $p$ и выполнено \eqref{eq10*} тогда и только тогда, когда все $\Omega_{p,\tilde r}^{X}$ $-$ геодезические пространства неотрицательной по Александрову кривизны.

\end{corollary}

\noindent{\bf 4. Характеризация $CAT(0)$ предкасательных пространств через неравенство Птолемея и выпуклость по Буземанну.} Практически одновременно с работой \cite{BN} была опубликована статья T. Foertsch, A. Lytchak, V. Schroeder \cite{FLS}, в которой CAT(0)-пространства были охарактеризованы как птолемеевы геодезические пространства вы\-пуклые по Буземанну.

Напомним, что метрическое пространство $(X,d)$ называется \emph{птолемеевым}, если следующее \emph{неравенство Птолемея}
\begin{equation}\label{eq13}
d(x,y)d(u,v)\le d(x,u)d(y,v)+d(x,v)d(y,u)
\end{equation}
имеет место для всех $x,y,u,v\in X.$

\begin{definition}\emph{\cite [с. 187]{Pa}}
Геодезическое пространство $(X,d)$ называется выпуклым по Буземанну (пространством Буземанна), если для любых двух аффинно параметризованных геодезических $\gamma: [a,b]\rightarrow X$ и $\gamma': [a', b']\rightarrow X$ отображение $D_{\gamma,\gamma'}:[a,b]\times [a',b']\rightarrow\mathbb R,$ определенное как $$D_{\gamma,\gamma'}(t,t')=d(\gamma(t),\gamma(t'))$$ является выпуклым.
\end{definition}

Мы не даем полное описание всех терминов, входящих в определение 4, отсылая читателя к монографии \cite{Pa}.

\begin{theorem}\emph{\cite{FLS}}
Метрическое пространство $X$ является СAT(0)-пространством тогда и только тогда, когда $X$ птолемеево и выпукло по Буземанну.
\end{theorem}

Целью настоящего раздела является построение инфинитезимального аналога теоремы 7.

\begin{remark}
В работе \cite{FLS} геодезическое пространство $(X,d)$ называется выпуклым по Буземанну, если для любых двух аффинно параметризованных геодезических $\beta: [a,b]\rightarrow X$ и $\gamma: [a,b]\rightarrow X$ отображение $t\mapsto d(\beta(t), \gamma(t))$ является выпуклым. Используя утверждение 8.12 из \cite{Pa}, легко показать, что такое определение эквивалентно определению 4.
\end{remark}

\textbf{Птолемеевость предкасательных.}
Условие птолемеевости предкасательных пространств было найдено в \cite{BDPtol}. Его легко получить, используя неравенство Птолемея \eqref{eq13} и приведенный выше принцип переноса. Аналогично функциям $A_1$ и $A_2$ зададим для пространства $(X,d,p)$ функцию $A_3$ как
\begin{equation}\label{eq14}
A_3 (w,x,y,z)=\frac{d(x,w)d(y,z)+d(x,z)d(y,w)-d(x,y)d(w,z)}{(\delta_{p}(w,x,y,z))^{2}}
\end{equation}
при $(w,x,y,z)\ne (p,p,p,p)$ и $A_{3}(p,p,p,p)=0.$

\begin{lemma}\emph{\cite{BDPtol}}
Пусть $(X,d,p)$ $-$ метрическое пространство с отмеченной точкой $p.$ Любое предкасательное пространство $\Omega_{p,\tilde r}^{X}$ является птолемеевым тогда и только тогда, когда \begin{equation}\label{eq15}\liminf_{w,x,y,z\to p}A_{3}(w,x,y,z)\ge 0.
\end{equation}
\end{lemma}

\textbf{Выпуклость предкасательных пространств по Буземанну.} Вопрос о выпуклости по Буземанну предкасательных пространств к общим метрическим пространствам ранее не исследовался и, видимо, не может быть получен как простое следствие принципа переноса. Для того, чтобы определить ``инфинитезимальную выпуклость по Буземанну'' мы будем использовать характеризацию выпуклости по Буземанну в терминах срединных точек.

Напомним, что если $\gamma: [a,b]\rightarrow X$ есть геодезическая,  $x=\gamma(a), y=\gamma(b),$ то образ отрезка $[a,b]$ при отображении $\gamma$ называется \emph{геодезическим сегментом} и обозначается $[x,y].$ Таким образом, геодезические сегменты в $X$ это в точности подмножества $X,$ изометричные отрезкам прямой.

\begin{lemma}
Пусть $(X,d)$ $-$ метрическое пространство, $a,b,c\in X$ и $[a,b], [b,c]$ $-$ геодезические сегменты, лежащие в $X.$ Если имеет место равенство
\begin{equation}\label{eq26}
d(a,b)+d(b,c)=d(a,c),
\end{equation}
то множество $[a,b] \cup [b,c]$ $-$ геодезический сегмент в $X.$ В частности, если $d(a,b)=$ $=d(b,c)=\frac{1}{2}d(a,c),$ то точка $b$ принадлежит геодезическому сегменту $[a,b] \cup [b,c]$ и является срединной для точек $a$ и $c$.
\end{lemma}

Мы не будем приводить формальное доказательство этой леммы. Напомним только, что (с точностью до параметризации) геодезическая в метрическом пространстве есть кривая, длина которой совпадает с расстоянием между её концами.

\begin{lemma}
Пусть $(X,d)$ $-$ геодезическое пространство. Следующие утверждения эк\-вивалентны.
\newline \emph{(i)} $X$ выпукло по Буземанну;
\newline \emph{(ii)} Для любых трех точек $x^{0}, x^{1}, y\in X$ выполнено неравенство
\begin{equation}\label{eq16} d(m^{x},y)\le\frac{1}{2}(d(x^{0},y)+d(x^{1},y)),
\end{equation} где $m^{x}$  $-$ срединная точка для $x^{0},x^{1}.$
\end{lemma}
\begin{proof}
Как известно, $X$ является выпуклым по Буземанну тогда и только тогда, когда для любой точки $y\in X$ и любого геодезического сегмента $[x^0, x^1],$ лежащего в $X,$ выполнено \eqref{eq16} (см. утверждение 8.12 из \cite{Pa}). Следовательно $(ii)\Rightarrow (i)$ доказано. Для проверки $(i)\Rightarrow (ii)$ допустим, что $X$ $-$ выпукло по Буземанну, $x^0, x^1, m_x, y\in X$ и $m_x$ $-$ срединная точка для $x^0, x^1.$ Если $x^0=x^1,$ то $x^0=m_x=x^1$ и \eqref{eq16} очевидно. Пусть $x^0\ne x^1.$ Так как $X$ $-$ геодезическое, то существуют геодезические сегменты $[x^0, m_x]$ и $[m_x, x^1],$ а так как $m_x$ $-$ срединная точка для $x^0, x^1,$ то $d(x^0, m_x)=d(x^1, m_x)=\frac{1}{2}d(x^0, x^1).$ Отсюда, по лемме 7, следует, что множество $[x^0, m_x]\cup [m_x, x^1]$ $-$ геодезический сегмент и мы опять можем использовать утверждение 8.12 из \cite{Pa}.
\end{proof}

Пусть теперь $(X,d,p)$ $-$ метрическое пространство срединно выпуклое в отмеченной точке $p$ в смысле определения 3. Пусть $\tilde x^{0}=\{x_{n}^{0}\}_{n\in\mathbb N}$ и $\tilde x^{1}=\{x_{n}^{1}\}_{n\in\mathbb N}$ две последовательности, принадлежащие $\tilde X$ и такие, что $$\lim_{n\to\infty}x_n^{0} =\lim_{n\to\infty}x_n^{1}=p. $$ Через $\tilde m^{x}=\{m_{n}^{x}\}_{n\in\mathbb N}$ будем обозначать инфинитезимальную срединную точку для $\tilde x^{0},$ $\tilde x^{1}.$

\begin{definition}
 Пусть $(X,d,p)$ $-$ метрическое пространство с отмеченной точкой $p.$ Будем говорить, что пространство $(X,d,p)$ выпукло по Буземанну в точке $p,$ если оно срединно выпукло в $p$ и для любых $\tilde x^{0}, \tilde x^{1}, \tilde y\in\tilde X,$ сходящихся к $p,$ и любой $\tilde m^{x}=\{m_{n}^{x}\}_{n\in\mathbb N}$  имеет место соотношение
\begin{equation}\label{eq17}
(d(m_{n}^{x}, y)-\frac{1}{2}(d(x_{n}^{0},y_{n})+d(x_{n}^{1},y_{n})))_{+}=o(\delta_{p}(x_{n}^{0},y_{n},x_{n}^{1})),
\end{equation} где
\begin{equation}\label{eq18}
(t)_{+}=\frac{|t|+t}{2}, \, t\in\mathbb R.
\end{equation}
\end{definition}

Следующая теорема дает достаточное условие выпуклости по Буземанну сепарабельных касательных пространств.

\begin{theorem}
Пусть $(X,d,p)$ $-$ метрическое пространство с отмеченной точкой $p.$ Если $X$ выпукло по Буземанну в точке $p,$ то любое сепарабельное касательное пространство $\Omega_{p,\tilde r}^{X}$ является геодезическим пространством выпуклым по Буземанну.
\end{theorem}

\begin{proof}
Пусть  $X$ выпукло по Буземанну в точке $p,$  $\Omega_{p,\tilde r}^{X}$ $-$ сепарабельное касательное пространство с метрикой $\rho$, а $\tilde X_{p,\tilde r}$ $-$ максимальное самостабильное семейство, соответствующее $\Omega_{p,\tilde r}^{X}.$ Из определения 5 следует, что $X$ $-$ срединно выпукло в точке. Следовательно, по лемме 3, $\Omega_{p,\tilde r}^{X}$ является геодезическим. Нужно показать, что $\Omega_{p,\tilde r}^{X}$ выпукло по Буземанну. Пусть $\gamma^{i}, \beta\in \Omega_{p,\tilde r}^{X}, i=0,1$ и пусть $\mu^{\gamma}$ $-$ срединная точка для $\gamma^0, \gamma^1$. В соответствии с леммой 8, выпуклость по Буземанну пространства $\Omega_{p,\tilde r}^{X}$ будет доказана, если
\begin{equation}\label{eq23}
\rho(\mu^{\gamma}, \beta)\le \frac{1}{2}(\rho(\gamma^{0}, \beta)+\rho(\gamma^{1}, \beta)).
\end{equation}
Если $\gamma^0 = \gamma^1,$ то $\mu^{\gamma}=\gamma^0=\gamma^1$ и \eqref{eq23} очевидно. Пусть $\gamma^0 \ne \gamma^1.$ Рассмотрим $\tilde x^{i}~=~\{x_n^{i}\}_{n\in\mathbb N}, \tilde y=\{y_n\}_{n\in\mathbb N}, \tilde z^{\gamma}=\{z_n^{\gamma}\}_{n\in\mathbb N},$ принадлежащие $\tilde X_{p,\tilde r}$ такие, что $\pi(\tilde x^{i})=\gamma^{i}, \pi(\tilde y)=\beta, \pi(\tilde z^{\gamma})=\mu^{\gamma}, i=0,1,$
где $\pi$ $-$ отображение проектирования $\tilde X_{p,\tilde r}$ на $\Omega_{p,\tilde r}^{X}$ (см. \eqref{diag}). Проверим, что $\tilde z^{\gamma}$ $-$ инфинитезимальная срединная точка для $\tilde x^{0}$ и $\tilde x^{1}.$ Положим $\alpha=\pi(\tilde p),$ где $\tilde p$ $-$ стационарная последовательность $(p,p,p,...).$ Используя определение метрики $\rho,$ находим
\begin{equation*}
\lim_{n\to\infty}\frac{\left|d(x_{n}^{1}, z_{n}^{\gamma})-\frac{1}{2}d(x_{n}^{1},x_{n}^{0})\right|}{d(x_{n}^{1},p)\vee d(x_{n}^{0},p)}=\lim_{n\to\infty}\frac{\frac{\left|d(x_{n}^{1}, z_{n}^{\gamma})-\frac{1}{2}d(x_{n}^{1},x_{n}^{0})\right|}{r_n}}{\frac{d(x_{n}^{1}, p)\vee d(x_{n}^{0},p)}{r_n}}=
\end{equation*}
\begin{equation*}
=\frac{1}{\rho(\alpha,\gamma^{1})\vee \rho(\alpha,\gamma^{0})}\left|\rho(\gamma^{1},\mu^{\gamma})-\frac{1}{2}\rho(\gamma^{1},\gamma^{0})\right|.
\end{equation*}
Так как $\gamma^{0}\ne\gamma^{1},$ то $0<\rho(\gamma, \gamma^{1})\vee \rho(\gamma, \gamma^{0})<\infty,$ а так как $\mu^{\gamma}$ $-$ срединная точка для $\gamma^{0}$ и $\gamma^{1},$ то
\begin{equation*}
\rho(\gamma^{1}, \mu^{\gamma})=\frac{1}{2}\rho(\gamma^{1},\gamma^{0}).
\end{equation*}
Следовательно
\begin{equation}\label{eq24}
\lim_{n\to\infty}\frac{\left|d(x_{n}^{1}, z_{n}^{\gamma})-\frac{1}{2}d(x_{n}^{1},x_{n}^{0})\right|}{d(x_{n}^{1},p)\vee d(x_{n}^{0},p)}=0.
\end{equation}

Аналогично находим
\begin{equation}\label{eq25}
\lim_{n\to\infty}\frac{\left|d(x_{n}^{0}, z_{n}^{\gamma})-\frac{1}{2}d(x_{n}^{1},x_{n}^{0})\right|}{d(x_{n}^{1},p)\vee d(x_{n}^{0},p)}=0.
\end{equation}
Равенства \eqref{eq24}, \eqref{eq25} и \eqref{eq9} показывают, что $\tilde z^{\gamma}$ $-$ инфинитезимальная срединная точка для $\tilde x^{0}$ и $\tilde x^{1}.$  Воспользовавшись \eqref{eq17}, получаем
\begin{equation*}
0=\lim_{n\to\infty}\frac{(d(z_n^{\gamma}, y_n)-\frac{1}{2}(d(x_n^{0},y_n)+d(x_n^{1},y_{n})))_{+}}{d(x_n^{0},p)\vee d(x_n^{1},p)\vee d(y_n,p)}=
\end{equation*}

\begin{equation*}
=\frac{(\rho(\mu^{\gamma},\beta)-\frac{1}{2}(\rho(\gamma^{0},\beta)+\rho(\gamma^{1},\beta)))_{+}}{\rho(\gamma^{0},\alpha)\vee \rho(\gamma^{1},\alpha)\vee \rho(\beta,\alpha)}.
\end{equation*}

Так как знаменатель последней дроби конечное положительное число, то
\begin{equation*}
(\rho(\mu^{\gamma},\beta)-\frac{1}{2}(\rho(\gamma^{0},\beta)+\rho(\gamma^{1},\beta)))_{+}=0,
\end{equation*}
что равносильно \eqref{eq23}.
\end{proof}

Теоремы 7, 8 и лемма 6 дают следующее

\begin{corollary}
Пусть $(X,d,p)$ $-$ метрическое пространство с отмеченной точкой $p$. Предположим, что любое предкасательное пространство $\Omega_{p,\tilde r}^{X}$ является сепарабельным и касательным. Тогда, если выполнено неравенство \eqref{eq15} и $X$ выпукло по Буземанну в точке $p,$ то каждое $\Omega_{p,\tilde r}^{X}$ является пространством класса CAT(0).
\end{corollary}

\begin{theorem}
Пусть $(X,d,p)$ $-$ метрическое пространство с отмеченной точкой $p.$ Если все предкасательные пространства $\Omega_{p,\tilde r}^{X}$ являются геодезическими, выпуклыми по Буземанну пространствами, то $X$ выпукло по Буземанну в точке $p.$
\end{theorem}

\begin{proof}
Пусть все $\Omega_{p,\tilde r}^{X}$ являются геодезическими, выпуклыми по Буземанну пространствами, но $X$ не выпукло по Буземанну в точке $p.$ Так как все $\Omega_{p,\tilde r}^{X}$ $-$ геодезические, то по лемме 4, пространство $X$ $-$ срединно выпукло в точке $p.$ Так как $X$ не является выпуклым по Буземанну в точке $p,$ то найдутся сходящиеся к $p$ последовательности $\tilde x^{0}, \tilde x^{1}$ с инфинитезимальной срединной точкой $\tilde m^{x},$ сходящаяся к $p$ последовательность $\tilde y$  и постоянная $c>0$ такие, что
\begin{equation}\label{eq27}
(d(m_n^{x},y_n)-\frac{1}{2}(d(x_n^{0}, y_n)+d(x_n^{1}, y_n)))_{+}>c\delta_{p}(x_n^{0}, y_n, x_n^{1})
\end{equation}
при всех $n\in\mathbb N.$ Покажем, что правая часть в \eqref{eq27} положительна при всех достаточно больших $n.$ Действительно, если это не так, то переходя к подпоследовательности можем считать, что $x_n^{0}=x_n^{1}=y_n=p$ для всех $n,$ но тогда $m_n^{x}=p$ для всех достаточно больших $n$ (см. замечание 1). Подставляя $p$ вместо всех переменных в \eqref{eq27}, получаем ложное неравенство $(0)_{+}>c\cdot 0.$ Таким образом, можно разделить правую и левую часть в \eqref{eq27} на $\delta_{p}(x_n^{0}, y_n, x_n^{1})$ и переходя, если нужно, к подпоследовательности считать, что
\begin{equation}\label{eq28}
0<c<\lim_{n\to\infty}\frac{d(m_n^{x}, y_n)-\frac{1}{2}(d(x_n^{0}, y_n)+d(x_n^{1}, y_n))}{\delta_{p}(x_n^{0}, y_n, x_n^{1})}.
\end{equation}
Из \eqref{eq9} следует неравенство
\begin{equation*}
\limsup_{n\to\infty}\frac{d(m_n^{x}, p)}{\delta_{p}(x_n^{0}, y_n, x_n^1)}<\infty.
\end{equation*}
Значит, переходя еще раз, если необходимо, к подпоследовательности можем считать, что $\tilde p, \, \tilde x^{1}=\{x_n^{1}\}_{n\in\mathbb N}, \tilde y=\{y_n\}_{n\in\mathbb N}, \tilde x^{0}=\{x_n^{0}\}_{n\in\mathbb N}, \tilde m^{x}=\{m_{n}^{x}\}_{n\in\mathbb N}$ $-$ взаимностабильные относительно нормирующей последовательности $\tilde r=\{r_n\}_{n\in\mathbb N}$ с $r_n=$ $=\delta_{p}(x_n^{0}, y_n, x_n^{1}), n\in\mathbb N.$ Пусть $\tilde X_{p,\tilde r}$ $-$ максимальное самостабильное семейство такое, что $\{\tilde x^{1}, \tilde x^{0}, \tilde y, \tilde m^{x}\}\subseteq \tilde X_{p,\tilde r}.$ Обозначим через $\Omega_{p,\tilde r}^{X}$ соответствующее семейству $\tilde X_{p,\tilde r}$ предкасательное пространство. Обозначим $$\beta^{0}=\pi(\tilde x^{0}), \beta^{1}=\pi(\tilde x^{1}), \gamma=\pi(\tilde y), \mu^{\beta}=\pi(\tilde m^{x}), \alpha=\pi(\tilde p).$$ Тогда \eqref{eq28} можно записать как

\begin{equation}\label{eq29}
0<c<\rho(\mu^{\beta}, \gamma)-\frac{1}{2}(d(\beta^{0},\gamma)+d(\beta^{1},\gamma)).
\end{equation}

В силу леммы 8, из последнего неравенства следует, что $\Omega_{p,\tilde r}^{X}$ не выпукло по Буземанну. Полученное противоречие завершает доказательство теоремы. \end{proof}

Теорема7, теорема 9 и лемма 6 дают следующее

\begin{corollary}
Пусть $(X,d,p)$ $-$ метрическое пространство с отмеченной точкой $p$. Если все предкасательные пространства $\Omega_{p,\tilde r}^{X}$ являются СAT(0)-прост\-ранствами, то $X$ выпукло по Буземанну в точке $p$ и выполняется неравенство \eqref{eq15}.
\end{corollary}

\begin{theorem}
Пусть $(X,d,p)$ $-$ метрическое пространство с отмеченной точкой $p$. Предположим, что любое $\Omega_{p,\tilde r}^{X}$ является сепарабельным и касательным. Тогда следующие утверждения эквивалентны.
\newline (i) Каждое $\Omega_{p,\tilde r}^{X}$ является пространством класса CAT(0).
\newline(ii) Пространство $X$ выпукло по Буземанну в точке $p$ и $\mathop{\liminf}\limits_{w,x,y,z\to p}A_{3}(w,x,y,z)\ge 0,$ где функция $A_3$ определена в \eqref{eq14}.
\end{theorem}

В доказанных выше теоремах 3, 8 и 10 фигурирует свойство сепарабельности касательных пространств $\Omega_{p,\tilde r}^{X}.$ В связи с этим возникает вопрос о свойствах пространства $X,$ гарантирующих такую сепарабельность. Одним из таких свойств является так называемое \emph{свойство удвоения}. Напомним, что метрическое пространство $X$ обладает \emph{свойством удвоения} (\emph{doubling property}), если существует $C \ge 1$ такое, что для любого $A\subseteq X$ с $\textrm{diam} A<\infty$ найдется не более чем $C$ множеств $B\subseteq X,$ образующих покрытие $A$ и удовлетворяющих условию $$\textrm{diam }B\le \frac{1}{2}\textrm{diam} A.$$ Используя утверждение 4 из \cite{ADK} можно доказать следующее

\bigskip

\noindent \textbf{Утверждение 1.} \emph{Пусть $X$ $-$ метрическое пространство, обладающее свойством удвоения. Тогда для любой точки $p\in X$ любое предкасательное пространство $\Omega_{p,\tilde r}^{X}$ является сепарабельным.}

Отсюда следует, что сепарабельными являются предкасательные к любому подмножеству конечномерного нормированного пространства. Интересно отметить, что в работе Д. Дордовского  \cite[доказательство теоремы 2.3]{Dor} было построено подмножество $X$ гильбертова пространства $l_2$ суммируемых с квадратом последовательностей, имеющее по крайней мере мере одно не сепарабельное предкасательное пространство. Таким образом, пространство предкасательное к сепарабельному метрическому пространству не обязательно сепарабельно.

\medskip

\textbf{Viktoriia Bilet}

Institute of Applied Mathematics and Mechanics of NASU, R. Luxemburg str. 74, Donetsk 83114, Ukraine

\textbf{E-mail:} biletvictoriya@mail.ru

\bigskip

\textbf{Oleksiy Dovgoshey}

Institute of Applied Mathematics and Mechanics of NASU, R. Luxemburg str. 74, Donetsk 83114, Ukraine

\textbf{E-mail:} aleksdov@mail.ru

\begin{thebibliography}{9}

\bibitem{ADK} F. Abdullayev, O. Dovgoshey, M. K\"{u}\c{c}\"{u}kaslan, \emph{Compactness and boundedness of tangent spaces to metric spaces}, Beitr. Algebra Geom., \textbf{51} (2010), № 2, 547 -- 576.

\bibitem {BN} I. Berg, I. Nikolaev, \emph{Quasilinearization and curvature of Aleksandrov spaces}, Geom. Dedicata, \textbf{133} (2008), 195 -- 218.

\bibitem {Bilet} V. Bilet, \emph{Geodesic tangent spaces to metric spaces}, Ukr. Mat. J., \textbf{64} (2012), № 9, 1273 -- 1281. (in Russian)

\bibitem {BD} V. Bilet, O. Dovgoshey, \emph{Isometric embeddings of pretangent spaces in $E^{n}$}, preprint available at http://arxiv.org/abs/1205.2335.

\bibitem {BDPtol} V. Bilet, O. Dovgoshey, \emph{Metric betweenness, Ptolemaic spaces and isometric embeddings of pretangent spaces in $R$}, J. Math. Sci., New York, \textbf{182} (2012),  № 4, 22 -- 36.

\bibitem {BBI} D. Burago, Yu. Burago, S. Ivanov, \emph{A course in metric geometry},
American Mathematical Society, Providence, RI, (2001).

\bibitem {Dor} D. Dordovskyi, \emph{Metric tangent spaces to Euclidean spaces }, J. Math. Sci., New York, \textbf{179} (2011),  № 2, 229 -- 244.

\bibitem {DM} O. Dovgoshey, O. Martio, \emph{Tangent spaces to metric spaces}, Reports in Math., Helsinki Univ., \textbf{480} (2008), 20 p.

\bibitem {DM1} O. Dovgoshey, O. Martio, \emph{Tangent spaces to general metric spaces}, Rev. Roumaine Math. Pures. Appl., \textbf{56} (2011), № 2,  137 -- 155.

\bibitem {FLS} T. Foertsch, A. Lytchak, V. Schroeder,  \emph{Nonpositive curvature and the Ptolemy inequality}, Int. Math. Res. Not., (2007), № 22, 15 p.

\bibitem {LP} N. Lebedeva, A. Petrunin, \emph{Curvature bounded below: a definition a la Berg-Nikolaev}, Electronic research announcements in math. sci., \textbf{17} (2010),  122 -- 124.

\bibitem {Pa} A. Papadopoulos, \emph{Metric spaces, convexity and nonpositive curvature},
European Mathematical Society, (2005).

\bibitem {Sato} T. Sato, \emph{An alternative proof of Berg and Nikilaev's characterization of CAT(0)-spaces via quadrialateral inequality}, Arch. Math., \textbf{93} (2009), 487 -- 490.


\end{thebibliography}
\end{document}